\documentclass[11pt]{amsart}
\usepackage{amsmath,amssymb}
\usepackage[all,cmtip]{xy}

\begin{document}

\newtheorem{thm}{Theorem}[section]
\newtheorem{lem}[thm]{Lemma}
\newtheorem{prop}[thm]{Proposition}
\newtheorem{cor}[thm]{Corollary}
\newtheorem{defn}[thm]{Definition}
\newtheorem*{remark}{Remark}

\numberwithin{equation}{section}

\newcommand{\Z}{{\mathbb Z}} 
\newcommand{\Q}{{\mathbb Q}}
\newcommand{\R}{{\mathbb R}}
\newcommand{\C}{{\mathbb C}}
\newcommand{\N}{{\mathbb N}}
\newcommand{\FF}{{\mathbb F}}
\newcommand{\fe}{\overline{\mathbb F}}
\newcommand{\fq}{\mathbb{F}_q}
\newcommand{\feq}{\overline{\mathbb F}_q}

\newcommand{\rmk}[1]{\footnote{{\bf Comment:} #1}}

\renewcommand{\mod}{\;\operatorname{mod}}
\newcommand{\ord}{\operatorname{ord}}
\newcommand{\TT}{\mathbb{T}}
\renewcommand{\i}{{\mathrm{i}}}
\renewcommand{\d}{{\mathrm{d}}}
\newcommand{\HH}{\mathbb H}
\newcommand{\Vol}{\operatorname{vol}}
\newcommand{\area}{\operatorname{area}}
\newcommand{\tr}{\operatorname{tr}}
\newcommand{\norm}{\mathcal N} 
\newcommand{\intinf}{\int_{-\infty}^\infty}
\newcommand{\ave}[1]{\left\langle#1\right\rangle} 
\newcommand{\Var}{\operatorname{Var}}
\newcommand{\Prob}{\operatorname{Prob}}
\newcommand{\sym}{\operatorname{Sym}}
\newcommand{\disc}{\operatorname{disc}}
\newcommand{\CA}{{\mathcal C}_A}
\newcommand{\cond}{\operatorname{cond}} 
\newcommand{\lcm}{\operatorname{lcm}}
\newcommand{\Kl}{\operatorname{Kl}} 
\newcommand{\leg}[2]{\left( \frac{#1}{#2} \right)}  

\newcommand{\sumstar}{\sideset \and^{*} \to \sum}

\newcommand{\LL}{\mathcal L} 
\newcommand{\sumf}{\sum^\flat}
\newcommand{\Hgev}{\mathcal H_{2g+2,q}}
\newcommand{\USp}{\operatorname{USp}}
\newcommand{\conv}{*}
\newcommand{\dist} {\operatorname{dist}}
\newcommand{\CF}{c_0} 
\newcommand{\kerp}{\mathcal K}

\newcommand{\fs}{\mathfrak S}
\newcommand{\rest}{\operatorname{Res}} 
\newcommand{\af}{\mathbb A} 
\newcommand{\Ht}{\operatorname{Ht}}
\title
[Square-free values of polynomials] {Square-free values of
polynomials over the rational function field }
\author{ Ze\'ev Rudnick}
\address{Raymond and Beverly Sackler School of Mathematical Sciences,
Tel Aviv University, Tel Aviv 69978, Israel}
\email{rudnick@post.tau.ac.il}

\begin{abstract}
We study representation of square-free polynomials in the polynomial
ring $\fq[t]$ over a finite field $\fq$ by polynomials in
$\fq[t][x]$. This is a function field version of the well studied
problem of representing square-free integers by integer polynomials,
where it is conjectured that   a separable polynomial  $f\in \Z[x]$
takes infinitely many square-free values, barring some simple
exceptional cases, in fact  that the integers $a$ for which $f(a)$
is square-free have a positive density. We show that if $f(x)\in
\fq[t][x]$ is separable, with square-free content,
 of bounded degree and height, then as $q\to \infty$,  for almost all monic polynomials $a(t)$, the polynomial $f(a)$ is square-free.
\end{abstract}

\thanks{We thank Lior Rosenzweig for his comments.
 This work was conceived during the ERC Research Period on
Diophantine Geometry in Pisa during September 2012. The research
leading to these results has received funding from the European
Research Council under the European Union's Seventh Framework
Programme (FP7/2007-2013) / ERC grant agreement n$^{\text{o}}$
320755.}
\date{\today}

\maketitle
\section{Introduction}
Let $\fq$ be a finite field of $q$ elements. We wish to study
representation of square-free polynomials in the polynomial ring
$\fq[t]$ by polynomials in $\fq[t][x]$. This is a function field
version of the well studied problem of representing square-free
integers by integer polynomials, where it is conjectured that   a
separable polynomial (that is, without repeated roots) $f\in \Z[x]$
takes infinitely many square-free values, barring some simple
exceptional cases, in fact  that the integers $a$ for which $f(a)$
is square-free have a positive density.
The  problem is most difficult when $f$ is irreducible. The
quadratic case was solved by Ricci  \cite{Ricci}. For cubics,
Erd\"{o}s \cite{Erdos} showed that there are infinitely many
square-free values, and Hooley \cite{Hooley 1967} gave the result
about positive density. Beyond that nothing seems known
unconditionally for irreducible $f$, for instance it is still not
known that $a^4 +2$ is infinitely often square-free. Granville
\cite{Granville} showed that the ABC conjecture completely settles
this problem.
An easier problem which has recently been solved is to ask how often
an irreducible polynomial $f\in \Z[x]$ of degree $d$ attains values
which are free of $(d-1)$-th powers, either when evaluated at
integers or at primes, see \cite{Erdos, Hooley77, Nair1, Nair2,
Helfgott, Browning, HB, Reuss}.


In this note we study  a function field version of this problem.
Given a   polynomial $f(x) = \sum_j \gamma_j(t)x^j\in \FF[t][x]$
which is separable, that is with no repeated roots in any extension
of $\fq(t)$,
we want to know how often is $f(a)$ square-free in $\fq[t]$ as $a$
runs over (monic) polynomials in $\fq[t]$.

We want to rule out polynomials like $f(x,t)=t^2x$ for which
$f(a(t),t)$ can never be square-free. To do so, recall that the
content $c\in \fq[t]$ of a polynomial $f\in \fq[t][x]$ as above is
defined as the greatest common divisor of the coefficients of $f$:
$c = \gcd(\gamma_0,\dots, \gamma_\ell)$. A polynomial is {\em
primitive} if $c=1$, and any $f\in \fq[t][x]$ can be written as $f=c
f_0$ where $f_0$ is primitive. If the content $c$ is not square-free
then $f(a)$ can never be square-free.

For any field $\FF$, let
\begin{equation}
M_n(\FF)=\{a\in \FF[t]: \deg a=n,a \mbox{ monic} \} \;,
\end{equation}
so that $\#M_n(\fq)=q^n$. Defining
\begin{equation}
\mathcal S_f(n)(\FF) = \{a\in M_n(\FF): f(a) \mbox{ is square-free }
\} \;,
\end{equation}
we want to study the frequency
\begin{equation}
\frac {\#\mathcal S_f(n)(\fq)}{\#M_n(\fq)}
\end{equation}
in an appropriate limit.

There are two possible limits to take: Large degree ($n\to \infty$)
while keeping the constant field $\fq$ fixed, or large constant
field ($q\to \infty$) while keeping $n$ fixed. The large degree
limit ($q$ fixed, $n\to \infty$) was investigated by Ramsay
\cite{Ramsay} who showed:
 \begin{thm}\label{Ramsay thm}
Assume $f\in \fq[t][x]$ is separable and  irreducible.
Then
\begin{equation*}
\frac{\#\mathcal S_f(n)(\fq)}{\#M_n(\fq)} = c_f + O_{f,q}(\frac
{1}n),\quad \mbox{as } n\to \infty \;,
\end{equation*}
with
$$c_f=\prod_P (1-\frac{\rho_f(P^2)}{|P|^2})\;,
$$
the product over prime polynomials $P$, and for any polynomial $D\in
\fq[t]$, $\rho_f(D)=\#\{C\bmod D: f(C)=0\bmod D\}$. The implied
constant depends on $f$ and on the finite field size $q$. The
density $c_f$ is positive if and only if there is some $a\in \fq[t]$
such that $f(a)$ is square-free.
\end{thm}
Ramsay actually counts all polynomials up to degree $n$, and does
not impose the monic condition. See also Poonen \cite{Poonen Duke}
for multi-variable versions.

Ramsay's theorem is proved by an elementary sieve argument, with one
crucial novel ingredient due to Elkies
to deal with the contribution of large primes to the sieve, which is
completely unavailable in the number field case; in Granville's work
\cite{Granville}, the ABC conjecture plays an analogous r\^ole.

In this note we deal with the large finite field limit, of $q\to
\infty$ while $n$ is fixed. Here it makes little sense to fix the
polynomial $f$, so we also allow variable $f$, as long as restrict
the degree (in $x$) and height, where for a polynomial $f(x,t) =
\sum_j \gamma_j(t)x^j\in \FF[t][x]$, the height is
$\Ht(f)=\max_j\deg \gamma_j(t)$.

We will show
\begin{thm}\label{thm finite field}
 For all separable $f\in \fq[t][x]$ with square-free content,   as $q\to \infty$,
\begin{equation}
\frac{\#\mathcal S_f(n)(\fq)}{ \#M_n(\fq)} = 1 +O(\frac {(n\deg
f+\Ht(f))\deg f}{q})\;,
\end{equation}
the implied constant absolute.
\end{thm}
Thus if we fix $n$, the degree and the height, as $q\to \infty$ for
almost all $a\in M_n(\fq)$ the polynomials $f(a)$ are square-free.
For instance, the number of  $a(t)\in M_n(\fq)$ for which $a(t)^4+2$
is square-free is, for $q$ odd, $q^n+O(nq^{n-1})$.

Note that since primes (irreducibles)  have positive density among
all monic polynomials of given degree in $\fq[t]$, we in particular
find that for almost all primes $P\in \fq[t]$ of given degree, the
polynomial $f(P)$ is square-free as $q\to \infty$.


Remark: It is possible to have primitive, separable $f$ with no
square-free values, for instance take
\begin{equation}
f(x) = \prod_{\alpha,\beta\in \fq} (x-\alpha t-\beta)=x^{2q}+\dots
\;.
\end{equation}
Then for all $a\in \fq[t]$, $f(a)$ is divisible by
$(\prod_{\gamma\in \fq}(t-\gamma))^2=(t^q-t)^2$. Indeed, if we fix
$\gamma\in \fq$, any $a\in \fq[t]$ is congruent modulo
$(t-\gamma)^2$ to some $\alpha t+\beta$ and hence $f(a) \equiv
f(\alpha t+\beta) = 0 \bmod (t-\gamma)^2$. Thus we need to impose
some restriction on the degree of $f$ in Theorem~\ref{thm finite
field}.

 Theorem~\ref{thm finite field} is a consequence
of a purely algebraic result, valid over {\em any} field $\FF$.

\begin{thm}\label{algebraic thm}
Suppose $f\in \FF[t][x]$ is  separable over $\FF(t)$ and has
square-free content. Then $\mathcal S_f(n)(\FF)$ is the complement
of a \underline{proper} Zariski-closed hypersurface of the affine
$n$-dimensional space $M_n(\FF)$, of degree $D\leq 2(n\deg f + \Ht
f) \deg f$.
\end{thm}


  Theorem~\ref{algebraic thm}  implies that the number of
$a\in M_n(\fq)$ for which $f(a)$ is not square-free is at most
$Dq^{n-1}$, where $D$ is the total degree of an equation defining
the hypersurface. Indeed, if $h \in \fq[X_1,\dots, X_m]$ is a
non-zero polynomial of total degree at most $D$,  then the number of
zeros of $h(X_1,\dots, X_m)$ in $\fq^m$ is at most $Dq^{m-1}$. This
is an elementary fact, seen by fixing all variables but one (cf
\cite[\S 4, Lemma 3.1]{Schmidt}). Hence Theorem~\ref{thm finite
field} follows.

\section{Proof of Theorem~\ref{algebraic thm}}

\subsection{The primitive case}

 We write
\begin{equation}\label{form of f}
f(x,t) = \gamma_0(t) + \gamma_1(t) x + \dots +\gamma_\ell(t) x^\ell
\end{equation}
with $\gamma_j(t)\in \FF[t]$, and $\gamma_\ell(t)\neq 0$. We first
assume that $f(x,t)$ is primitive, that is $\gcd(\gamma_j(t)) = 1$.
Denote by
\begin{equation}
\Delta_f(t)=\disc_x f(x,t)
\end{equation}
 the discriminant of $f(x)$ as a polynomial of degree $\ell$ with coefficients in $\FF[t]$; it is a universal polynomial
with integer coefficients in $\gamma_0(t), \dots, \gamma_\ell(t)$:
\begin{equation}
\Delta_f(t) = \mbox{Poly}_\Z(\gamma_0(t),\dots ,\gamma_\ell(t))\in
\FF[t] \;.
\end{equation}
Separability of $f$  (over $\FF(t)$) is equivalent to the
discriminant not being the zero polynomial: $\Delta_f(t)\neq 0$.

 The key observation is that
$f(a)\in \FF[t]$ being square-free is equivalent to requiring that
the polynomial $ t\mapsto f(a(t),t)$ does not have any multiple
zeros (in any extension of the field $\FF$). This is in fact a
polynomial condition, that is a polynomial system of equations for
the coefficients $a_0,a_1,\dots, a_{n-1}$ of $a(t)=a_0+a_1t+\dots +
a_{n-1}t^{n-1}+t^n$ which is given by the vanishing of the
discriminant:
\begin{equation}\label{Hypersurface}
\disc f(a(t),t)  =0 \;.
\end{equation}
It suffices to show that this equation defines a {\em proper}
hypersurface. 

 Before doing so, we bound the degree $D$ of the hypersurface
 \eqref{Hypersurface}: For $f(x,t)$ as in \eqref{form of f},
$f(a(t),t)$ is a polynomial in $t$ of degree
 \begin{equation}
 \deg   f(a(t),t)\leq n\deg f+\max \deg \gamma_j = n\deg f +
 \Ht(f) \;.
 \end{equation}
The coefficients are polynomials in the $a_j$ of degree at most
$\deg f$.
 Now the discriminant of a polynomial $\sum_{j=0}^m h_j t^j$ is
 homogeneous in the coefficients $h_j$ of degree $2m-2$. Hence
 $a\mapsto \disc f(a(t),t)=\sum_k \delta_k \prod a_i^{k_i}$ has
 total degree at most
 \begin{equation}\label{degree D}
 D\leq 2(n\deg f+\Ht(f))\deg f \;.
 \end{equation}
It remains to show that the equation \eqref{Hypersurface} is
nontrivial.

The condition that the polynomial $f(a(t))$ has multiple zeros  is
that there is some $\rho\in \fe$ (an algebraic closure of $\FF$)
with
\begin{equation}\label{common zeros 1}
f(a(\rho),\rho)=0, \quad \frac{\partial f}{\partial
x}(a(\rho),\rho)\cdot  a'(\rho)  +\frac{\partial f}{\partial
t}(a(\rho),\rho) =0 \;.
\end{equation}

 we define
\begin{equation}
W=\{(\rho,\vec a)\in \af^1\times \af^{n }: \eqref{common zeros 1}
\mbox{
 holds}\} \;.
\end{equation}
We have a fibration of $W$ over the $\rho$ line $\af^1$ and a map
$\phi: W\to \af^{n }$, the restriction of the projection
$\af^1\times \af^{n }\to \af^{n }$,
\begin{equation}
\xymatrix{ &{W\subset \af^1\times \af^{n }} \ar[dl]_{\pi}
\ar[dr]^{\phi}
&\\
{\af^1}  && {\af^{n }} }
\end{equation}
and the solutions of \eqref{common zeros 1} are precisely $\phi(W)$.

We will show that generically the fiber $\pi^{-1}(\rho)$ has
dimension $n-2$ and for at most finitely many $\rho$ the dimension
is $n-1$. Therefore we obtain  that $\dim W=n-1$. Since the
solutions of \eqref{common zeros 1} are precisely $ \phi(W)$, it
follows that $\dim \phi(W)\leq n-1$. This will conclude the proof of
Theorem~\ref{algebraic thm} in the primitive case.

We note that for primitive polynomials,  $f(x,\rho) = \sum_j
\gamma_j(\rho)x^j$ is not the zero polynomial for any $\rho\in \fe$.
Thus for each $\rho \in \fe$, the condition $f(a(\rho),\rho)=0$
constrains $a$ to solve an equation $a(\rho)=\beta$, where $\beta\in
\fe$ is on of the at most $\deg f$ roots of $f(x,\rho)$.

We separate into two cases:    The singular case when
$\frac{\partial f}{\partial x}(a(\rho),\rho) =0$
 and the generic case when we require $\frac{\partial f}{\partial x}(a(\rho),\rho) \neq
 0$

The singular case 
implies that $\beta$ is a multiple zero of the polynomial
$f(x,\rho)$, that is that  $\rho$ is a zero of the discriminant
$\Delta_f(t)$, which is not identically zero (since we assume $f$ is
separable) and hence there are only finitely many possibilities for
such $\rho$. Given one of those $\rho$,  then we need $a(t)$ to
satisfy $a(\rho)=\beta$, i.e.
\begin{equation}
a_0+a_1\rho+\dots +a_{n-1}\rho^{n-1} +\rho^n=\beta
\end{equation}
which is a (non-degenerate) linear equation, and therefore carves
out an $n-1$-dimensional subspace of $a$'s. Thus the singular locus
consists of at most finitely many hyperplanes, and hence if
non-empty has dimension $n-1$.

In the generic case, we substitute $a(\rho)=\beta$ into
\eqref{common zeros 1}  to get a system
\begin{equation}
a(\rho) = \beta, \quad a'(\rho) = -\frac{\frac{\partial f}{\partial
x}(\beta,\rho)}{\frac{\partial f}{\partial t}(\beta ,\rho)}
\end{equation}
that is
\begin{equation}
\begin{split}
a_0 +&a_1\rho+a_2\rho^2+\dots +a_{n-1}\rho^{n-1}  =-\rho^n +\beta \\
\, &a_1 + a_2\cdot 2\rho +\dots + a_{n-1} \cdot (n-1)\rho^{n-2} =
-n\rho^{n-1} - \frac{\frac{\partial f}{\partial
x}(\beta,\rho)}{\frac{\partial f}{\partial t}(\beta ,\rho)}
\end{split}
\end{equation}
which is clearly of rank $2$. Hence the fibers $\pi^{-1}(\rho)$ have
dimension $n-2$.

\subsection{The general case}
We now relax the primitivity condition. Write $f(x,t) =
c(t)f_0(x,t)$ where $f_0(x,t)=\sum_j \gamma_j^{(0)}(t)x^j$  is
primitive, and $c(t)\in \fq[t]$ is square-free. Since $c(t)$ is
square-free, we obtain that $f(a(t),t)$ is square-free if and only
if $f_0(a(t),t)$ is square-free and coprime to $c(t)$. Now
$f_0(a(t),t)$ being square-free is the condition $\disc
f_0(a(t),t)\neq 0$. For $f_0(a)$ to not be coprime to $c$ is the
algebraic condition on vanishing of the resultant
\begin{equation}
R=\rest(c(t),f_0(a(t),t)) \;.
\end{equation}
 Thus the set of $a\in M_n$ so that $f(a)$ is square-free
 is the complement of the hypersurface
\begin{equation}\label{disc R}
\disc f_0(a(t),t) \cdot R(t) = 0 \;.
\end{equation}
We wish to show that this is a non-zero equation and to bound its
total degree.

We have established above that the discriminant equation $\disc
f_0(a(t),t) =0$ is nontrivial,   of total degree
\begin{equation}
D_0\leq 2(n\deg f_0+\Ht(f_0))\deg f_0 =2(n\deg f+\Ht(f_0))\deg f
\end{equation}
in $a_0,\dots,a_n$.

We wish to show that the resultant $R$ is not identically zero.
Assuming (as we may) that $c(t)$ is monic, we can write the
resultant as a product over the zeros of $c(t)$
\begin{equation}\label{resultant expansion}
R
= 
\prod_{c(\alpha) = 0} f_0(a(\alpha),\alpha) = 
\prod_{c(\alpha) = 0} \sum_{j=0}^\ell \gamma_j^{(0)}(\alpha)
(a_0+\dots +   \alpha^n)^j \;.
\end{equation}
For each zero $\alpha$ of $c(t)$, let $  \ell(\alpha)=\deg
f_0(x,\alpha)$  be the degree of the polynomial $f_0(x,\alpha)\in
\feq [x]$,  which is not the zero polynomial by primitivity of
$f_0$. Then the total degree of $R$ is
\begin{equation}
L:= \sum_{c(\alpha)=0} \ell(\alpha) \leq \deg c \cdot \deg f
\end{equation}
and the coefficient of $a_0^L$ is $\prod_\alpha
\gamma_{\ell(\alpha)}^{(0)}(\alpha)$ which is non-zero. Hence $R$ is
non-zero and of degree $L$.

%

Finally, we compute the total degree of the equation \eqref{disc R}
is the sum of $D_0$ and $\deg R$, which is at most
\begin{equation}
 2(n\deg f+\Ht(f_0))\deg f +\deg c \cdot \deg f \leq 2(n\deg
 f+\Ht(f))\deg f
\end{equation}
since $\Ht(f) = \Ht(f_0) + \deg c$. This concludes the proof of
Theorem~\ref{algebraic thm}.

\end{document}